\documentclass{article}
\usepackage{amsmath}
\usepackage{amsthm}
\usepackage{amssymb,epsfig}
\usepackage[cp 1250]{inputenc}
\usepackage{graphics}
\usepackage{color}
\newtheorem{lem}{Lemma}[section]
\newtheorem{thm}[lem]{Theorem}
\newtheorem{definition}[lem]{Definition}
\newtheorem{corollary}[lem]{Corollary}
\newtheorem{lemma}[lem]{Lemma}
\newtheorem{remark}[lem]{Remark}

\newtheorem{proposition}[lem]{Proposition}

\newtheorem{observation}[lem]{Observation}

\title{Repetitions in  beta-integers}

\author{L{\!'}. Balkov\'a,  K. Klouda, E. Pelantov\'a
\footnote{e-mail: l.balkova@centrum.cz, karel@kloudak.eu, edita.pelantova@fjfi.cvut.cz}\\
\emph{Doppler Institute and Department of Mathematics FNSPE, }\\
\emph{Czech Technical University in Prague,}\\ \emph{ Trojanova
13, 120 00 Praha 2, Czech Republic}\\}
\date{}

\begin{document}

\maketitle

\begin{abstract}
Classical crystals are solid materials containing arbitrarily long
periodic repetitions of a~single motif. In this Letter, we study
the maximal possible repetition of the same motif occurring in
$\beta$-integers -- one dimensional models of quasicrystals.
We are interested in
$\beta$-integers realizing only a~finite number of distinct
distances between neighboring elements. In~such~a~case, the
problem may be reformulated in terms of combinatorics on words as
a~study of the index of infinite words coding $\beta$-integers. We
will solve a~particular case for $\beta$ being a~quadratic
non-simple Parry number.
\end{abstract}
\section{Introduction}

\medskip
This Letter takes up the study of $\beta$-integers initiated by
the investigation of their asymptotic properties in~\cite{BaGaPe}.
Similarly as in the previous Letter, we restrict our consideration
to $\beta$-integers realizing only a~finite number of distinct
distances between neighbors; $\beta$ is then called a~Parry
number. For Parry numbers,  the set of $\beta$-integers forms
a~discrete aperiodic Delone set with a~self-similarity factor
$\beta$ and of finite local complexity. It follows herefrom that
$\beta$-integers are suitable for modeling materials with aperiodic
long range order, the so-called quasicrystals~\cite{FrGaKr}. Classical
crystals are solid materials containing arbitrarily long periodic
repetitions of a~single motif. Quasicrystals do not share this property.

In this Letter, we are interested in the maximal
possible repetition of one motif occurring in $\beta$-integers.
It turns out to be suitable to reformulate and study this problem
in terms of combinatorics on words.

For Parry numbers, coding distinct distances between
neighboring nonnegative $\beta$-integers with distinct letters,
one obtains a~right-sided infinite word $u_\beta$ over a~finite
alphabet. The reformulation of our task in the language of
combinatorics on words has the following reading: For a~given
factor $w$ of the infinite word $u=u_\beta$, find the longest
prefix $v$ of the infinite periodic word $w^\omega=wwwwww \ldots$
such that $v$ occurs as a~factor in~$u=u_\beta$. The ratio of the
lengths of $v$ and $w$ is called the index of the factor $w$ in
$u=u_\beta$ and is denoted by $\text{ind}(w)$. Let us note that
$\text{ind}(w)$ is not necessarily an integer. Denote by $k$ the
lower integer part ${{\lfloor \text{ind}(w)\rfloor}}$ of the index
of $w$, then the word $w^k$, i.e., the concatenation of $k$ words
$w$, is usually called the maximal integer power of $w$.

The index of any infinite word $u$ can be  naturally defined as
$$\text{ind}(u)=\sup\{\text{ind}(w)\bigm | w \ \text{factor of} \
u\}.$$ Explicit values of the index are known only for few classes
of infinite words. The index of Sturmian words has been studied
in many papers (\cite{Be1}, \cite{JuPi}, \cite{MiPi}, \cite{MaPe}),
the complete solution to the problem was provided
independently by Carpi and de Luca~\cite{CaLu} and by Damanik and
Lenz~\cite{DaLe}. Recently, the index of infinite words has
reinforced its importance: Damanik in~\cite{Da} considers
discrete one-dimensional Schr\"{o}dinger operators with aperiodic
potentials generated by primitive morphisms and he establishes
purely singular continuous spectrum with probability one provided
that the potentials (infinite words) have the index greater than
three. Let us stress that infinite words $u_\beta$
associated with Parry numbers belong to the class of
infinite words generated by primitive morphisms, too.

Here we study the index of infinite words $u_\beta$ for quadratic
non-simple Parry numbers $\beta$. These words are determined by
integer parameters $p,q$, where $ p>q\geq 1$. We provide an
explicit formula for $\text{ind}(u_\beta)$. In the particular case
of $p=q+1$, the infinite word $u_\beta$ is Sturmian and our result may be deduced also
from the well-known formula for the index of Sturmian words.
We have chosen the word $u_\beta$ associated with quadratic non-simple Parry
numbers $\beta$ for our study of the index of non-Sturmian words since for such infinite words,
we dispose of detailed knowledge on arithmetical properties of $\beta$-integers and combinatorial
properties of the associated infinite words $u_\beta$ (\cite{BaMa}, \cite{BaPeTu}).

The Letter is organized in the following way. In Section
\ref{kombinatorika}, we introduce necessary notions from
combinatorics on words and we cite a~relevant result on the index of
Sturmian words. In Section~\ref{betaintegers}, we provide the
background on infinite words $u_\beta$ coding $\beta$-integers for
$\beta$ being a~non-simple quadratic Parry number. In Section
\ref{integerpowers}, we determine the maximal integer power
occurring in $u_\beta$ (Theorem \ref{IntPower}).
Section~\ref{indexUbeta} is devoted to the index of $u_\beta$
(Theorem \ref{ind}) and to the comparison of our result with the
formula for the index of Sturmian words.

\section{Preliminaries}\label{kombinatorika}
An {\em alphabet} $\mathcal A$ is a~finite set of symbols, called {\em letters}. Throughout this paper, the binary alphabet ${\mathcal A}=\{0,1\}$ is used.
The string $w=w_1w_2
\ldots w_k$, where $w_i\in {\cal A}$ for each $i=1,2,\ldots,k$, is called a~{\em word of length} $k$
on ${\mathcal A}$. The length of $w$ is then $k$ and it is denoted by $|w|=k$. The set of all finite words together with the operation of concatenation
forms a~monoid; its neutral element is the empty word $\varepsilon$.
We denote this monoid ${\mathcal A}^*$.
An infinite sequence $u=u_0u_1u_2\ldots $ of symbols from the alphabet ${\mathcal A}$ is called an~{\em infinite word}. A~finite word $w$ is said to be a~{\em factor} of the (finite or infinite) word $v$ if
there exists a~finite word $v'$ and a~finite or infinite word $v''$ such that $v=v'wv''$.
If $v'$ is the empty word, then $w$ is called a~{\em prefix} of $v$, if $v''$ is the empty word,
then $w$ is a~{\em suffix} of $v$. If $v=v'w$, then $vw^{-1}=v'$, i.e., the word $vw^{-1}$ is obtained from $v$ by erasing its suffix $w$. The set of all factors of an infinite word $u$ is said to be the {\em language}
of $u$ and is denoted $\mathcal{L}(u)$. An infinite word $u$ is called
{\em recurrent} if every of its factors occurs
infinitely many times in $u$ and $u$ is called {\em uniformly recurrent} if for every of its factors $w$, the set of all factors in $\mathcal{L}(u)$ that do not contain $w$ as their factor is finite. In other words, every sufficiently long element of $\mathcal{L}(u)$ contains $w$ as its factor.

The number of factors of the infinite word $u$ gives us insight
 into its variability. The function $ \mathcal{C}_u : \mathbb{N}\mapsto \mathbb{N}$ that to every $n$ associates the number of distinct factors of length $n$ occurring in $u$ is called the {\em factor complexity} of the infinite word $u$.
An infinite periodic word $u=www\ldots$, where $w$ is a~finite word, is usually denoted $w^\omega$. Its factor complexity $\mathcal{C}_u$ is bounded; it is readily seen that $\mathcal{C}_u(n) \leq |w|$. Similarly, the factor complexity of an eventually
periodic word $u=w'w^\omega$, where $w', w$ are finite words, is bounded. The necessary and sufficient condition for an infinite word to be aperiodic is validity of the equation $\mathcal{C}_u (n) \geq n+1$ for all $n\in \mathbb{N}$. Infinite aperiodic words satisfying $\mathcal{C}_u (n)=n+1$ for all $n \in \mathbb N$ are called {\em Sturmian words}; these are thus infinite aperiodic words of the lowest possible factor complexity. Sturmian words are the best known aperiodic words; a~survey on their properties may be consulted in~\cite{Be}. In particular, any Sturmian word is uniformly recurrent.

For determination of the factor complexity of an infinite word $u$, an essential role is played by special factors. We recall that a~factor $w$ of an infinite word $u$ over a~binary alphabet $\{0,1\}$ is called {\em left special}
if $0w$ and $1w$ are both factors of $u$, $w$ is called {\em right
special} if $w0$ and $w1$ are both factors of $u$, and $w$ is said
to be {\em bispecial} if $w$ is both left special and right special.

\bigskip

The crucial notion of our study is the index of a~factor $w$ in a~given infinite word $u$.
Let us define first the integer power of $w$.
For any $k\in\mathbb{N}$, the $k$-{th \em power} of $w$ is the
concatenation of $k$ words $w$, usually denoted $w^k$.
Analogously for $r \in \mathbb Q, r\geq 1$,  we call a~word $v$
the $r$-{th \em power} of the word $w$ if there exists a~proper
prefix $w'$ of $w$ such that
$$v=\underbrace{w\dots w}_{\lfloor r \rfloor - times}
w'\quad \text{and} \quad r=\lfloor r \rfloor
+\frac{|w'|}{|w|}\,.$$
The $r$-th power of $w$ is denoted $w^r$.

Our aim is to find, for a~given factor $w$ of an infinite word
$u$, the highest power of $w$ occurring in~$u$.
We will be interested exclusively in aperiodic uniformly recurrent words.
Such words contain for every factor $w$ only $r$-th powers of $w$ with bounded $r$,
see~\cite{Lo}. Therefore, for aperiodic uniformly recurrent words $u$, it makes sense to define the index of $w$ in $u$ as
$$ \text{ind}(w)=\max \{ r \in \mathbb{Q} \mid w^r \in
\mathcal{L}(u)\}\,.$$ The word $w^{\text{ind}( w)}$ is called the {\em maximal
power} of $w$ in $u$, the word $w^{\lfloor \text{ind} ( w)\rfloor}$
the {\em maximal integer power} of $w$ in $u$. The index of the infinite word $u$ is defined as
$$\text{ind}(u)=\sup\{\text{ind}(w)\bigm | w \in \mathcal{L}(u)\} \
.$$ Let us remark that an aperiodic uniformly recurrent word $u$
can have an infinite index; even among Sturmian words, one can find words with an infinite index.
The language of every Sturmian word is characterized by an irrational parameter $\alpha \in (0,1)$,
called {\em slope}. If $\alpha$ is the slope of a~Sturmian word $u$, then the word obtained by exchanging letters in $u$ has the slope $1-\alpha$ and has evidently the same index as $u$. Consequently, we may assume without loss of generality that $\alpha>\tfrac{1}{2}$. In order to determine the index of a~Sturmian word, we need to express $\alpha$ in the form of its continued fraction. Since
$\tfrac{1}{2}<\alpha<1$, the continued fraction of $\alpha$ equals $[0;1, a_2,a_3, \ldots ]$. The results of \cite{CaLu} and~\cite{DaLe} say that the index of a~Sturmian word $u$ with slope
$\alpha$ is equal to
\begin{equation}\label{sturm}
\text{ind} (u)=2+ \sup \Bigl\{ a_{n+2} + \frac{q_n-2}{q_{n+1}}\mid
n \in \mathbb{N} \Bigr\},
\end{equation}
where $q_n$ is the denominator of the $n$-th convergent of $\alpha$.

Now, let us describe a~large class of uniformly recurrent words: fixed points of primitive morphisms. This class includes infinite words $u_\beta$ associated with Parry numbers $\beta$. The map
$\varphi : {\cal A}^* \mapsto {\cal A}^*$ is called a~{\em morphism} if
$\varphi(wv)=\varphi(w)\varphi(v)$ for every $w,v \in{\cal A}^*$.
One may associate with $\varphi$ the {\em morphism matrix} $M_\varphi$
satisfying
$$(M_\varphi)_{ab}=\text{number of letters $b$ occurring in $\varphi(a)$},$$
for any pair of letters $a,b\in {\mathcal A}$.

Knowing for any $a\in {\mathcal A}$ the number of letters $a$
occurring in a~factor $w$, we may obtain the same information for
$\varphi(w)$ by a~simple formula. We mention the formula only for the binary
alphabet ${\mathcal A}=\{0,1\}$ we are interested in. It
follows straightforwardly from the definition of the morphism matrix that for every
factor $w \in {\mathcal A}^{*}$
\begin{equation}\label{subst_matrix}
(|\varphi(w)|_0, |\varphi(w)|_1)=(|w|_0, |w|_1)M_\varphi.
\end{equation}
where $|v|_a$ denotes the number of letters $a$ occurring in a word
$v$. Clearly, a~similar formula holds for any finite alphabet.

A~morphism is said to be {\em primitive} if a~power of
$M_\varphi$ has all elements strictly positive. In other words, matrices of primitive morphisms fulfill the assumptions of the Perron-Frobenius theorem.

The action of the morphism $\varphi$ may be naturally extended to an infinite word $u=u_0u_1u_2 \dots$
by the prescription
$$\varphi(
u_0u_1u_2 \ldots )=\varphi(u_0)\varphi(u_1)\varphi(u_2) \ldots
$$
An infinite word $u$ is called a~{\em fixed point} of $\varphi$ if
$\varphi(u)=u$. It is known that any fixed point of a~primitive
morphism is uniformly recurrent and that any left eigenvector
corresponding to the dominant eigenvalue of $M_\varphi$ is
proportional to the densities of letters in any fixed point of $\varphi$.

\section{Infinite words associated with $\beta$-integers}\label{betaintegers}

Here, we provide the description of infinite words $u_\beta$
associated with quadratic non-simple Parry numbers in terms of
fixed points of primitive morphisms. We keep the notation from our
precedent Letter~\cite{BaGaPe}, where the number theoretical
background on $\beta$-integers and associated infinite words
$u_\beta$ is at disposal.  Nevertheless, to make this Letter self-contained,
we will recall all notions needed
for understanding of our results. In this section
we also deduce some important properties of $u_\beta$, in
particular, a~transformation generating bispecial factors of
$u_\beta$. Bispecial factors turn out to be essential for our main
aim -- determination of the maximal integer powers of factors and determination of the
index of~$u_\beta$.

A non-simple quadratic Parry number $\beta$ is the larger root of
$\mbox{$x^2-(p+1)x+p-q$}$, where $p>q \geq 1$. Its R\'enyi expansion of unity is
$d_{\beta}(1)=pq^{\omega}$. The set of $\beta$
integers has two distances between neighbors: $\Delta_0=1$ and
$\Delta_1=\beta-p$.  Consequently, the infinite word $u_\beta$
coding distances between neighboring $\beta$-integers is binary.

As we have already said, the reader may find the notions of the R\'enyi expansion of unity,
$\beta$-integers, distances between neighbors in $\mathbb Z_\beta$ etc.
in our precedent Letter~\cite{BaGaPe}. However, in order to follow the ideas in the sequel,
it is sufficient to know that $u_\beta$ is the unique fixed
point of a~morphism canonically associated with parameters $p,q$
characterizing non-simple quadratic Parry numbers $\beta$.
Therefore we will use the result of~\cite{Fa} for an equivalent definition of
$u_\beta$.

\begin{definition} Let $\beta$ is the larger root of
$\mbox{$x^2-(p+1)x+p-q$}$, where $p>q \geq 1$.  The unique fixed
point of the morphism
\begin{equation} \label{subst2}
\varphi(0)=0^p1,\quad \varphi(1)=0^q1\,
\end{equation}
will be denoted $u_\beta$.
\end{definition}

The infinite word $u_\beta$ starts as follows
\begin{equation}\label{u:subst2}
u_\beta=\underbrace{\underbrace{0^p1 \ \dots \ 0^p1}_{p
\hbox{\scriptsize\ times}} \ 0^q1 \ \underbrace{0^p1 \ \dots \
0^p1}_{p \hbox{\scriptsize\ times}} \ 0^q1 \ \dots \
\underbrace{0^p1 \ \dots \ 0^p1}_{p \hbox{\scriptsize\ times}} \
0^q1}_{p \hbox{\scriptsize\ times}}\ \underbrace{0^p1 \ \dots \
0^p1}_{q \hbox{\scriptsize\ times}}0^q1 \ \dots
\end{equation}
The morphism matrix $M_\varphi$ is
 $\left ( \begin{smallmatrix}
 p  &  1 \\
 q  &  1
  \end{smallmatrix} \right )$
and $\varphi$ is thus obviously primitive. Computing the left
eigenvector of $M_\varphi$ corresponding to the dominant
eigenvalue $\beta$, we get the densities $1 - \tfrac{1}{\beta}$
and $\tfrac{1}{\beta}$ of letters $0$ and $1$, respectively.

\begin{remark}
In the paper~\cite{BaMa}, it is shown that the factor complexity $\mathcal{C}$ of $u_\beta$
satisfies
$$ \hbox{if}\ \  p>q+1,\ \ \ \hbox{then}\ \
 \{\mathcal{C}(n+1) - \mathcal{C}(n)\mid n \in \mathbb{N} \}=\{1,2\},$$
$$ \hbox{if}\ \  p=q+1,\ \ \ \hbox{then}\ \
 \{\mathcal{C}(n+1) - \mathcal{C}(n)\mid n \in \mathbb{N} \}=\{1\}.~~~$$ Therefore, $u_\beta$ is Sturmian if and only if \
$p=q+1$.
\end{remark}

First of all, some simple, but very important
properties of the morphism $\varphi$ are observed.
\begin{observation}
\label{zero_blocks} Let $10^k1$ be a~factor of $u_\beta$, then $k=p$
or $k=q$.
\end{observation}
\begin{observation}\label{ancestors} Let $v$ be any factor of $u_\beta$
containing at least one 1. Then there exists $k_1, \ 0 \leq k_1 \leq p$,
such that $0^{k_1}1$ is a~prefix of $v$ and there exists $k_2, \ 0 \leq k_2 \leq p$,
such that $10^{k_2}$ is a~suffix of $v$. The fact that $\varphi(0)$ and $\varphi(1)$ end in
$1$ and contain only one letter $1$ implies that there exists
a~unique word $w$ in $\{0,1\}^*$ satisfying
$v=0^{k_1}1\varphi(w)0^{k_2}$. Clearly, $w$ is a~factor of $u_\beta$.
\end{observation}

Of significant importance is the map $T: \{0,1\}^* \rightarrow
\{0,1\}^*$ defined by
\begin{equation}\label{T_def}
T(w)=0^q1 \varphi(w)0^q.
\end{equation}

The map $T$ helps to generate bispecial factors 
that play a~crucial role in the determination of the
index of factors. Therefore, the rest of this section is devoted to the description of properties of $T$.

\begin{lemma}\label{T}
Let $T$ be the map defined in~\eqref{T_def}.
\begin{enumerate}
\item \label{inLanguage} For every $w \in {\mathcal L}(u_\beta)$, it holds that $T(w) \in {\mathcal L}(u_\beta)$.
\item \label{Textensions} Let $w$ be a~factor of $u_\beta$ and let $a,b \in {\mathcal A}$, then
$awb \in {\mathcal L}(u_\beta)$ if and only if $aT(w)b \in {\mathcal
L}(u_\beta).$
\item \label{Tbispecials} Let $v$ be a~bispecial factor
of $u_\beta$ containing at least one letter 1, then there exists
a~unique factor $w$ such that $v=T(w)$.
\item \label{prefix} Let
$w,v$ be factors of $u_\beta$, then $w$ is a~prefix of $v$ if and
only if $T(w)$ is a~prefix of $T(v)$.
\item \label{suffix} Let
$w,v$ be factors of $u_\beta$, then $w$ is a~suffix of $v$ if and
only if $T(w)$ is a~suffix of $T(v)$.
\end{enumerate}
\end{lemma}
\begin{proof}
{\em 1.} Take an arbitrary factor $w \in {\mathcal L}(u_\beta)$.
Then $w$ is extendable to the right, and, since $u_\beta$ is
recurrent, $w$ is also extendable to the left. In other words, there
exists $a,b \in \{0,1\}$ such that $awb$ is also a~factor of
$u_\beta$. As $u_\beta$ is a~fixed point of $\varphi$, the image
$\varphi(awb)$ belongs to ${\mathcal L}(u_\beta)$. Finally, $T(w)$
is a~factor of $u_\beta$ because $T(w)$ is a~subword of
$\varphi(awb)$.

{\em 2.} Let $1w1$ be a~factor of ${\mathcal L}(u_\beta)$, then
$01w1$ is as well a~factor of $u_\beta$. Applying $\varphi$, we
learn that $\varphi(01w1)=\varphi(0)T(w)1$ is a~factor of $u_\beta$,
which proves that $1T(w)1$ belongs to ${\mathcal L}(u_\beta)$. The
other cases $0w0, 0w1, 1w0$ are analogous.

Let $0T(w)1 \in {\mathcal L}(u_\beta)$, using
Observation~\ref{zero_blocks}, the word $v=10^p1\varphi(w)0^q1$ is
also a~factor of $u_\beta$. Applying Observation~\ref{ancestors},
$v=1\varphi(0w1)$ and $0w1$ is an element of ${\mathcal L}(u_\beta)$.
All the other cases $0T(w)0, 1T(w)0,
1T(w)1$ are similar.

{\em 3.} Observation~\ref{zero_blocks} implies that each bispecial factor
$v$ containing at least one letter 1 has the prefix $0^q1$ and the suffix
$10^q$. According to Observation~\ref{ancestors}, there exists
a~unique $w$ such that $v=T(w)$.

{\em 4.} The implication $\Rightarrow$ is obvious noticing that
$0^q$ is a~prefix of $\varphi(a)$ for $a \in \{0,1\}$. The opposite
implication $\Leftarrow$ follows taking into account that $\varphi(1)$ is not a~prefix
of $\varphi(0)$ and $\varphi(0)$ is not a~prefix of $\varphi(1a)$ for any $a \in \{0,1\}$.

{\em 5.} The implication $\Rightarrow$ is obvious noticing that
$0^q1$ is a~suffix of $\varphi(a)$ for $a \in \{0,1\}$. The opposite
implication $\Leftarrow$ follows taking into account that $1\varphi(1)$ is not a~suffix of $\varphi(0)$ and $\varphi(0)$ is not a~suffix
of $\varphi(x1)$ for any $x \in \{0,1\}^{*}$.
\end{proof}

\section{Integer powers in $u_\beta$}\label{integerpowers}
Even if we want to describe the maximal integer powers of factors of $u_\beta$, it turns out to be useful to study first the relation between bispecial factors and the maximal rational powers of factors.
\begin{lemma}\label{BSandPOWER}  Let $u$ be an infinite uniformly
recurrent word over an alphabet $\mathcal A$. Let $w^kw'$ be its factor
for some proper prefix $w'$ of $w$ and some positive integer $k$.
Let us denote by {\bf P1}, \ {\bf P2}, \ {\bf P3} the following statements:
\begin{description}
\item[P1] The factor $w$ has the maximal index in $u$ among all factors of
  $u$ with the same length $|w|$ and $w^kw'$ is the maximal power of $w$ in $u$.
\item[P2] There exist $a,b \in \{0,1\}$ such that
$$aw^kw'b \in {\mathcal L}(u) \quad \text{and} \quad w'b \ \text{is not a~prefix of $w$} \quad \text{and} \quad  a \ \text{is not a~suffix of $w$}.$$
\item[P3] All the following factors are bispecial:
$$w', ww', www', \dots, w^{k-1}w'.$$
\end{description}
Then $\bf P1$ implies $\bf P2$ and $\bf P2$ implies $\bf P3$.

\end{lemma}
\begin{proof}
${\mathbf P1} \Rightarrow {\mathbf P2}:$ \quad  As $u$ is recurrent, there exists $a$ such that $aw^kw'$ is a~factor of $u$.
Since $w^kw'$ is the maximal power of $w$, the letter $a$ is not a~suffix of $w$,
otherwise the factor  $awa^{-1}$ (usually called a~conjugate of $w$)
 would have a~larger index than $w$.
On the other hand, if $w^kw'b$ is a~factor of $w$, then $w'b$ is not a~prefix of $w$,
otherwise it contradicts the fact that $w^kw'$ is the maximal power of $w$ in $u$.

${\mathbf P2} \Rightarrow {\mathbf P3}:$ \quad
Since $w'$ is a~proper prefix of $w$, there exists $x \in \mathcal A$ such that $w'x$ is a~prefix of $w$. Denote by $y$ the last letter of $w$. Obviously, $x \not=b$ and $y \not=a$. As $aw^jw'x$ is a~prefix of $aw^kw'b$ and $yw^jw'b$ is a~suffix of $aw^kw'b$, both $aw^jw'x$ and $yw^jw'b$ are in ${\mathcal L}(u)$ for all $j, \ 0 \leq j \leq k-1$. It follows that all factors listed in $\bf P3$ are bispecial.
\end{proof}

In this section, our aim is to describe the maximal integer powers occurring in $u_\beta$. Since the letter $0$ has the maximal index $p$, we may restrict our consideration to $k$-th powers of factors of $u_\beta$ with $k \geq p$.
Crucial for the determination of the index of $u_\beta$ are Propositions~\ref{essential} and~\ref{essential_one}.
\begin{proposition}\label{essential}
Let $p,q$ be integers, $p>q\geq 1$, and  $u_\beta$ be the fixed point of the morphism
 \eqref{subst2}. Assume $p>3$.
Let $w$ be a~factor of $u_\beta$ containing at least two $1$s and $w'$ be a~proper prefix of $w$.
Denote $v=w^kw'$ for some $k \in \mathbb N, \ k \geq p$. If there exist $a,b \in \{0,1\}$ so that
$$avb \in {\mathcal L}(u_\beta)\quad \text{and} \quad w'b \ \text{is not a~prefix of $w$} \quad \text{and} \quad
a \ \text{is not a~suffix
of $w$},$$
then there exist a~unique $\tilde{w}$ of length $\geq 2$
                    and a~proper prefix $\tilde{w}'$ of $\tilde{w}$ such that
                    \begin{equation}\label{Bletchley}
                        w=0^q1 \ \varphi(\tilde{w}) (0^q1)^{-1}
                        \quad\text{and}\quad v=T(\tilde{v})=T(\tilde{w}^k \tilde{w}');
                    \end{equation}
moreover,               $$\ a\tilde{v}b \in {\mathcal L}(u_\beta) \quad \text{and} \quad \tilde{w}'b \ \text{is not a~prefix of $\tilde{w}$} \quad \text{and} \quad a \ \text{is not a~suffix of $\tilde{w}$}.$$
\end{proposition}
In order to prove Proposition~\ref{essential}, we will use the following lemma.
\begin{lemma}\label{additional}
Let $p,q$ be integers, $p>q\geq 1$, and  $u_\beta$ be the fixed point of the morphism
 \eqref{subst2}.  The following statements
hold:
\begin{enumerate}

\item If \  $0(x1)^\ell x0  \in {\mathcal L}(u_\beta)$ \ for some integer $\ell\geq 2$, then
$\ell=2$ and $x=0^q$.
\item  If \   $1(x0)^\ell x1 \in {\mathcal L}(u_\beta)$
 \ for some integer $\ell\geq p-1$  and  $p\leq 2q$,
  then $x$ is the empty word $\varepsilon$.
\end{enumerate}
\end{lemma}
\begin{proof}
\begin{enumerate}
\item At first, we exclude the case when $x$ contains a non-zero letter.
 Suppose that the letter $1$ occurs in $x$. Since the factors $0x1$ and $1x0$ belong to ${\mathcal
 L}(u_\beta)$, it follows that $x$ is bispecial. By Lemma~\ref{T} Item~{\em \ref{Tbispecials}.}, $x$ starts in $0^q1$ and ends in $10^q$.
 Therefore $10^q10^q1 \in {\mathcal L}(u_\beta)$.  As
  $10^q10^q1=1 \varphi(11)$, we have according to
  Observation~\ref{ancestors} that $11\in  {\mathcal L}(u_\beta)$
  -- a contradiction.

   Now consider $x=0^s$ for some $s\in \mathbb{N}$. Then
   $0x1, 1x1 \in  {\mathcal L}(u_\beta)$, which implies by Observation~\ref{zero_blocks} that $s=q$. If $\ell$
   was at least $3$, then $1x1x1=10^q10^q1 \in {\mathcal
   L}(u_\beta)$, which leads to the same contradiction as before.

 \item  Again we start with the case of $x$ containing the letter $1$. 
 Since factors $1x0$ and $0x1$ belong to the language ${\mathcal L}(u_\beta)$, 
 it follows that $x$ is bispecial. Hence $x$ starts in $0^q1$ and ends in $10^q$. Since $10^q00^q1 \in {\mathcal L}(u_\beta)$, by Observation~\ref{zero_blocks}, we have $p=2q+1$.
This contradicts the assumption $p\leq 2q$.

Suppose now that $x=0^s$ for some  $s\in \mathbb{N}$. Since
$1(x0)^\ell x1=1(0^{s+1})^\ell 0^s1$ is a~factor of $u_\beta$,
Observation \ref{zero_blocks} gives that $(s+1)\ell +s \leq p$,
which is impossible if $s\geq 1$ and $\ell \geq p-1$.
Therefore $s=0$ and $x$ is the empty word.

\end{enumerate}
\end{proof}
\begin{proof}[Proof of Proposition~\ref{essential}]
The factor $w$ contains at least two $1$s. Since $w$ and $v=w^kw'$ satisfy Item {\bf P2} of Lemma~\ref{BSandPOWER}, both $ww'$ and $www'$ are bispecial, and therefore start in $0^q1$ and end in $10^q$. Consequently, their form is $ww'=T(x)$ and
$www'=T(y)$, where $x,y \not=\varepsilon$.
According to Lemma~\ref{T} Item~{\em \ref{suffix}.}, $x$ is a~suffix of $y$, i.e., $y=zx$ for some $z\not=\varepsilon$.
Observing $ww'=0^q1\varphi(x)0^q$ and $www'=0^q1\varphi(z)\varphi(x)0^q$, it follows directly that
$w=0^q1\varphi(z)(0^q1)^{-1}$.

Let us, at first, show that
\begin{enumerate}
\item
either $z$ is a~prefix of $x$,
\item or $z=x1$,
\item or $z=x0$.
\end{enumerate}
\noindent Assume $z=tdz'$ and $x=t(1-d)x'$ for a~word $t \in \{0,1\}^*$ and for a~letter $d$.
Then $w=0^q1\varphi(t)\varphi(d)\varphi(z')(0^q1)^{-1}$ is a~prefix of $ww'=0^q1\varphi(t)\varphi(1-d)\varphi(x')0^q$.
If $z'\not=\varepsilon$, we have a~contradiction immediately.
If $z'=\varepsilon$, then $t\not=\varepsilon$ knowing that $z$ contains at least two letters ($w$ contains at least two $1$s).
\begin{itemize}
\item If $d=1$, then $w=0^q1\varphi(t)$ and $ww'=0^q1\varphi(t)0^p1\varphi(x')0^q$, thus $w'$ starts in $0^p1$, which is not a~prefix of $w$ -- a~contradiction.
\item If $d=0$, then $w=0^q1\varphi(t)0^{p-q}$ and $ww'=0^q1\varphi(t)0^{p-q}0^{2q-p}1\varphi(x')0^q$, thus $w'$ starts in $0^{2q-p}1$, which is not a~prefix of $w$ because $p\not=q$ -- a~contradiction again.
\end{itemize}
\noindent The situation $z=xz'$ for some $z'$ of length $\geq 2$ cannot occur because it implies $|w|>|ww'|$.
Consequently, one of the situations $1., 2., 3.$ occurs.
\begin{enumerate}
\item If $z$ is a~prefix of $x$, i.e., $x=zx''$, then $ww'=0^q1\varphi(z)\varphi(x'')0^q$ and
$w=0^q1\varphi(z)(0^q1)^{-1}$, thus $w'=0^q1\varphi(x'')0^q$. Then
$v=w^{k}w'=0^q1\varphi(z^{k}x'')0^q$,
therefore $\widetilde{w}=z$, $\widetilde{w}'=x''$. As $w'$ is a~proper prefix of $w$, it follows by Lemma~\ref{T} Item~{\em \ref{prefix}} that $\widetilde{w}'$ is a~proper prefix of $\widetilde{w}$.
\item
Assume $z=x1$, then $w=0^q1\varphi(x)$ and $ww'=0^q1\varphi(x)0^q$, thus $w'=0^q$.
Then $v=0^q1\varphi((x1)^{k-1}x)0^q$. Since $w'0$ is not a~prefix of $w$ and $0$ is not a~suffix of $w$, the assumptions imply that $0v0 \in {\mathcal L}(u_\beta)$. By Observation~\ref{zero_blocks}, we have $1\varphi(0(x1)^{k-1}x0)\in {\mathcal L}(u_\beta)$. Since $k-1 \geq 3$, we deduce by Lemma~\ref{additional} that $0(x1)^{k-1}x0$ is not a~factor of $u_\beta$. It contradicts Observation~\ref{ancestors}. Hence, the case $z=x1$ does not occur.
\item If $z=x0$, then $w=0^q1\varphi(x)0^{p-q}$ and $ww'=0^q1\varphi(x)0^q$, thus $w'=0^{2q-p}$, which can happen only for $p\leq 2q$. Then $v=0^q1\varphi((x0)^{k-1}x)0^{q}$.
Since $w'1$ is not a~prefix of $w$ and $1$ is not a~suffix of $w$, the assumptions imply that $1v1 \in {\mathcal L}(u_\beta)$. Hence, $1\varphi(1(x0)^{k-1}x1) \in {\mathcal L}(u_\beta)$. However, by Lemma~\ref{additional}, it follows that $x=\varepsilon$. Then $z=0$, which contradicts the condition $|z|\geq 2$.
Thus, the case $z=x0$ does not occur.
\end{enumerate}
Let us finally note that the very last statement on the extensions of $\widetilde{v}$ follows from Lemma~\ref{T} Items~{\em \ref{Textensions}., \ref{prefix}.,} and~{\em \ref{suffix}.}
\end{proof}
\begin{remark}\label{exception}
Proposition~\ref{essential} does not take into account $u_\beta$ given by parameters $p=2,q=1$, $p=3,q=2$, and $p=3,q=1$. In the two first cases, $u_\beta$ is a~Sturmian word. Therefore, exclusion of the first two cases does not mean any loss.
For the case of $p=3, q=1$, in the proof of Proposition~\ref{essential}, we cannot exclude the situation $2.$; in this case,
Lemma~\ref{additional} Item 1. implies either the validity of~\eqref{Bletchley} or of
$$w=01\varphi(01)(01)^{-1} \quad v=T(01010).$$
\end{remark}

Proposition~\ref{essential} thus claims that for every factor $w \in {\mathcal L}(u_\beta)$ containing at least two $1$s
such that its $k$-th power $w^k$ is a~factor of $u_\beta$ with $k \geq p$, there exists a~shorter factor $\tilde{w}$ such that its $k$-th power $\tilde{w}$ is also in the language of $u_\beta$. As a~consequence, in order to determine the maximal integer power present in $u_\beta$, it is sufficient to study the index of factors $w$ containing only one letter $1$.

\begin{proposition}\label{essential_one}
Let $p,q$ be integers, $p>q\geq 1$, and  $u_\beta$ be the fixed point of the morphism
 \eqref{subst2}. Let $w$ be a~factor of $u_\beta$ containing one letter $1$
 and of the maximal index $\text{ind}(w)$ among all factors
 of length $|w|$ and such that $\text{ind}(w) \geq  p \geq 3$. Denote $k:=\lfloor \text{ind}(w)\rfloor$ and $v=w^kw'$ the maximal power of $w$.
 Then $$w=0^q1 \ \varphi(0) (0^q1)^{-1} \quad \text{and} \quad v=T(0^p)$$
 and $$\text{ind}(w)=p+\frac{2q+1}{p+1}.$$
\end{proposition}
\begin{proof}
According to Lemma~\ref{BSandPOWER}, $ww'$ is a~bispecial factor.
Lemma~\ref{T} Item {\em \ref{Tbispecials}.} claims that $ww'$ starts in $0^q1$. Therefore $w=0^q10^s$ with $s \in \{0, p-q\}$ (Observation~\ref{zero_blocks}). The case $s=0$ does not occur
since $www=0^q10^s0^q10^s0^q10^s \in {\mathcal L}(u_\beta)$ and $0^q10^q10^q1=0^q1\varphi(11)$, but $11 \not \in {\mathcal L}(u_\beta)$ -- a~contradiction to Observation~\ref{ancestors}. Hence $w=0^q10^{p-q}=0^q1\varphi(0)(0^q1)^{-1}$.
It remains to determine the form of $v$.
Again, since $ww'$ is bispecial, $ww'$ ends in $10^q$.
As $w'$ is a~prefix of $w$, at most one $1$ occurs in $w'$.
\begin{itemize}
\item
Suppose $w'$ contains $1$, then $w'$ starts in $0^q1$ and ends in $10^q$, thus $w'=0^q10^q$.
This is possible only in case when $p-q-1\geq q$, i.e., $p\geq 2q+1$.
Then, $v=(0^q10^{p-q})^{k}0^q10^q$. Consequently, $v=0^q1\varphi(0^{k})0^q=T(0^k)$.
On one hand, Observations~\ref{zero_blocks} and~\ref{ancestors} imply that $k \leq p$. On the other hand, since $v$ is the maximal power of $w$, it follows that $k\geq p$.
\item
Assume $w'$ does not contain $1$. Since $ww'$ ends in $10^q$ and $w=0^q10^{p-q}$, the only possibility for $w'$
is $w'=0^{2q-p}$. This comes in question only for $p\leq 2q$. Then $v=(0^q10^{p-q})^{k }0^{2q-p}=0^q1\varphi(0^{k-1})0^q=T(0^{k-1})$. The same arguments as in the previous case imply that $k-1=p$.
\end{itemize}
Clearly, $\text{ind}(w)=\frac{|v|}{|w|}=p+\frac{2q+1}{p+1}.$
\end{proof}
Let us state the main result of this section.
\begin{thm}\label{IntPower}
Let $p,q$ be integers, $p>q\geq 1$, and  $u_\beta$ be the fixed point of the morphism
 \eqref{subst2}. Assume $p \geq 3$.
\begin{itemize}
\item If $p\leq 2q$, then there exists a~factor $w \not=\varepsilon$ satisfying $w^{p+1} \in {\mathcal L}(u_\beta)$ and no $(p+2)$-nd power of any factor belongs to the language ${\mathcal L}(u_\beta)$.
\item If $p>2q$, then there exists a~factor $w \not=\varepsilon$ satisfying $w^{p} \in {\mathcal L}(u_\beta)$ and no $(p+1)$-st power of any factor belongs to the language ${\mathcal L}(u_\beta)$.
\end{itemize}
\end{thm}
\begin{proof}
Proposition~\ref{essential} implies that in order to determine the maximal integer power present in $u_\beta$, we can restrict our consideration to powers of factors containing only one letter $1$. When we compute the integer part
$\lfloor \text{ind}(w)\rfloor$ of such factors $w$ in Proposition~\ref{essential_one}, we find that the maximum is $p+1$ if $p\leq 2q$ and $p$ otherwise.
\end{proof}

\section{Index of $u_\beta$}\label{indexUbeta}
The task of this section is to compute the index of $u_\beta$, i.e.,
$$\text{ind}(u_\beta)=\sup\{\text{ind}(w)\bigm | w \in \mathcal{L}(u_\beta)\}.$$
We already know that $\text{ind}(u_\beta)\geq p$. Using Lemma~\ref{BSandPOWER}, it suffices to study rational powers $v=w^kw'$ of factors $w$ with the property $\bf P2$. As a~direct consequence of Propositions~\ref{essential} and~\ref{essential_one}, we have the following corollary.
\begin{corollary}
Let $p,q$ be integers, $p>q\geq 1$, and  $u_\beta$ be the fixed point of the morphism
 \eqref{subst2}. Assume $p>3$.
The index of $u_\beta$ is given by the following formula
$$
    \text{ind}(u_\beta)=\sup\{\text{ind}(w^{(n)})\bigm | n \in \mathbb N\},
$$
where
\begin{equation}\label{seq_of_w}
    w^{(0)}=0, \quad \quad w^{(n+1)}=0^q1\varphi(w^{(n)})(0^q1)^{-1}.
\end{equation}
Moreover, the maximal power of $w^{(n)}$ is $v^{(n)}$, where
\begin{equation}\label{seq_of_v}
    v^{(0)}=0^p, \quad \quad v^{(n+1)}=T(v^{(n)}).
\end{equation}
\end{corollary}
In the sequel, let us determine the index of $w^{(n)}$ for every $n \in \mathbb N$.
\begin{lemma}
The number of $0$s and $1$s in the words $w^{(n)}$ and $v^{(n)}$ satisfy
$$(|w^{(n)}|_0,|w^{(n)}|_1)=(1,0)M_\varphi^{n},$$
$$(|v^{(n)}|_0,|v^{(n)}|_1)=(p+1,\tfrac{2q+1-p}{q})M_\varphi^{n}-(1,\tfrac{2q+1-p}{q}),$$
where
 $M_\varphi=\left ( \begin{smallmatrix}
 p  &  1 \\
 q  &  1
  \end{smallmatrix} \right )$ is the morphism matrix.
\end{lemma}
\begin{proof}
As $w^{(n)}$ is a~conjugate of $\varphi(w^{(n-1)})$, the first formula holds by~\eqref{subst_matrix}.
Let us show the second one by induction on $n$.\\
\noindent For $n=0$, $$(|v^{(0)}|_0,|v^{(0)}|_1)=(p,0)=(p+1,\tfrac{2q+1-p}{q})-(1,\tfrac{2q+1-p}{q}).$$
\noindent For $n>0$, $$(|v^{(n)}|_0,|v^{(n)}|_1)=(|v^{(n-1)}|_0,|v^{(n-1)}|_1)M_\varphi+(2q,1)=$$
by the induction assumption,
$$=\left[(p+1,\tfrac{2q+1-p}{q})M_\varphi^{n-1}-(1,\tfrac{2q+1-p}{q})\right]M_\varphi+(2q,1)=(p+1,\tfrac{2q+1-p}{q})M_\varphi^{n}-(1,\tfrac{2q+1-p}{q}).$$
\end{proof}
Since the eigenvalues $\beta$ and $\beta'$ of $M_\varphi$ are roots of the Parry polynomial $x^2-(p+1)x+(p-q)$,
it is straightforward to show that $\vec{x_1}=(\beta-1,1)$ is a~left eigenvector of $M_\varphi$ corresponding to $\beta$
and $\vec{x_2}=(\beta'-1,1)$ is a~left eigenvector of $M_\varphi$ corresponding to $\beta'$.
The index of $w^{(n)}$ may be expressed as follows
$$\text{ind}(w^{(n)})=\frac{|v^{(n)}|}{|w^{(n)}|}=\frac{(p+1,0)M_\varphi^n\left(\begin{smallmatrix}1\\1\end{smallmatrix}\right)+
(0,\tfrac{2q+1-p}{q})M_\varphi^n\left(\begin{smallmatrix}1\\1\end{smallmatrix}\right)-(1,\tfrac{2q+1-p}{q})
\left(\begin{smallmatrix}1\\1\end{smallmatrix}\right)}{(1,0)M_\varphi^n\left(\begin{smallmatrix}1\\1\end{smallmatrix}\right)},$$
$$\text{ind}(w^{(n)})=p+1+\frac{\tfrac{2q+1-p}{q}(\alpha_1\vec{x_1}+\alpha_2\vec{x_2})M_\varphi^n\left(\begin{smallmatrix}
1\\1\end{smallmatrix}\right)-\tfrac{3q+1-p}{q}}{(\gamma_1 \vec{x_1}+\gamma_2\vec{x_2})M_\varphi^n\left(\begin{smallmatrix}1\\1\end{smallmatrix}\right)},$$
where $\alpha_1\vec{x_1}+\alpha_2\vec{x_2}=(0,1)$ and $\gamma_1\vec{x_1}+\gamma_2\vec{x_2}=(1,0)$.
Using the fact that $\vec{x_1}$ and $\vec{x_2}$ are eigenvectors of $M_\varphi$, we have
$$\text{ind}(w^{(n)})=p+1+\frac{\tfrac{2q+1-p}{q}\left(\alpha_1 \beta^n \vec{x_1}\left(\begin{smallmatrix}1\\1\end{smallmatrix}\right)+\alpha_2 \beta'^n\vec{x_2}\left(\begin{smallmatrix}1\\1\end{smallmatrix}\right)\right)-\tfrac{3q+1-p}{q}}{\gamma_1 \beta^n \vec{x_1}\left(\begin{smallmatrix}1\\1\end{smallmatrix}\right)+
\gamma_2\beta'^n\vec{x_2}\left(\begin{smallmatrix}1\\1\end{smallmatrix}\right)}.$$
It is easy to calculate that
$\alpha_1=\frac{1-\beta'}{\beta-\beta'}, \
\alpha_2=\frac{\beta-1}{\beta-\beta'}, \
\gamma_1=\frac{1}{\beta-\beta'}, \
\gamma_2=\frac{-1}{\beta-\beta'}$, and
$\vec{x_1}\left(\begin{smallmatrix}1\\1\end{smallmatrix}\right)=\beta$,
$\vec{x_2}\left(\begin{smallmatrix}1\\1\end{smallmatrix}\right)=\beta'$.
The final formula for the index of $w^{(n)}$ has the form
\begin{multline*}
    \text{ind}(w^{(n)})=p+1+\frac{\tfrac{2q+1-p}{q}\left((1-\beta')\beta^{n+1}-
    (1-\beta)\beta'^{n+1}\right)-\tfrac{3q+1-p}{q}(\beta-\beta')}{\beta^{n+1}-\beta'^{n+1}}=\\   
    =p+1+\frac{2q+1-p}{q}(1 - \beta') + \underbrace{\frac{\beta -
    \beta'}{q(\beta^{n+1} - {\beta'}^{n+1})}}_{>0}
    \underbrace{\left((2q + 1 - p){\beta'}^{n+1} - (3q +1 - p)
    \right)}_{\mathrm{A}(n)}.
\end{multline*}
Using the fact $0<\beta'<1<\beta$, we determine the limit
$$
    \lim_{n \to \infty}\text{ind}(w^{(n)})=p+1+\frac{2q+1-p}{q}(1-\beta')=p+1+\frac{2q+1-p}{\beta-1}.
$$
This limit is the supreme of $\{\text{ind}(w^{(n)})\bigm | n \geq 0\}$ if and only if $\mathrm{A}(n)<0$ for all $n \in \mathbb N$. It is an easy exercise to show that $\mathrm{A}(n)=(2q+1-p)({\beta'}^n - 1)-q<0$ for all $n \in \mathbb N$ if and
only if $p \leq 3q + 1$, otherwise
$\mathrm{A}(n)>0$ for all sufficiently large $n$.

Let us sum up the results in a~theorem.
\begin{thm}\label{ind} Let $p,q$ be integers, $p>q\geq 1$, and  $u_\beta$ be the fixed point of the morphism
 \eqref{subst2}. Assume $p>3$.
Then the index of $u_\beta$ satisfies

\noindent for $p \leq 3q + 1$ $$\text{ind}(u_\beta)=p+1+\frac{2q+1-p}{\beta-1},$$
otherwise there exists $n_0 \in \mathbb N$ such that $$\text{ind}(u_\beta)=\text{ind}(w^{(n_0)})>p+1+\frac{2q+1-p}{\beta-1}.$$
\end{thm}
\begin{remark}
Similarly as in the previous section, we have to treat the case of $p=3$ and $q=1$ separately.
According to Remark~\ref{exception}, we have to determine the index of $(w^{(n)})$ defined in~\eqref{seq_of_w}, but moreover the index of $(\hat{w}^{(n)})$ defined recursively by
$$\hat{w}^{(0)}=01\varphi(01)(01)^{-1}, \quad \quad \hat{w}^{(n+1)}=0^q1\varphi(\hat{w}^{(n)})(0^q1)^{-1}.$$
Using the same technique as before, we obtain
$$\sup\{\text{ind}(\hat{w}^{(n)})\bigm | n \in \mathbb N\} \quad = \quad \beta \quad  < \quad 4 \quad = \quad
\sup\{\text{ind}(w^{(n)})\bigm | n \in \mathbb N\}.$$
Hence, Theorem~\ref{ind} holds in fact also in this case.
\end{remark}

At the conclusion, let us compare in case of Sturmian words $u_\beta$ the formula for $\text{ind}(w^{(n)})$ with the formula~\eqref{sturm} for the index of general Sturmian words.
As we have already stated, $u_\beta$ is Sturmian if and only if
$p=q+1$, i.e., $\beta$ is the larger root of the polynomial $x^2-(p+1)x+1$.
For such $\beta$, we have
$$
    \text{ind}(u_\beta)=p+1+\frac{2p-1}{\beta-1}=\beta +1\,.
$$

In order to apply the formula from~\eqref{sturm}, we need to determine the slope $\alpha$ of the Sturmian word $u_\beta$. Since $\alpha >\frac{1}{2}$ is the density of the more frequent letter, according to Section~\ref{betaintegers}, $\alpha=1-\frac{1}{\beta}$. Let us use some basic properties of continued fractions available at any book on Number Theory to determine the continued fraction of this value. Since
$$ 1-\frac{1}{\beta}=\frac{1}{1+\frac{1}{\beta - 1}}=\frac{1}{1+\frac{1}{p-1 + 1-\frac{1}{\beta}}},$$
one obtains $\beta=[0;1,(p-1),1,(p-1),\ldots ]=[0;\overline{1,(p-1)} ]$.
Denominators $q_n$ of the convergents of $\beta$ fulfill therefore the following recurrent
relations
$$ q_{2n+1}=(p-1)q_{2n}+q_{2n-1} \qquad {\rm and } \qquad
q_{2n}=q_{2n-1}+q_{2n-2}$$ with initial values  $q_1=1,\ q_2=p,\
q_3=p+1$. By mathematical induction on $n$, it may be shown easily that
$$ q_{2n-1}=\frac{1}{\beta - \beta'}\Bigl(\beta^n -
\beta'^n\Bigr)\qquad {\rm and } \qquad q_{2n}=\frac{1}{\beta -
\beta'}\Bigl((1-\beta')\beta^{n+1} - (1-\beta)\beta'^{n+1}\Bigr).$$
As it holds for coefficients of the continued fraction of $\beta$ that $a_{2n-1}=1$ and $a_{2n}=p-1$, it suffices to consider even $n$ in~\eqref{sturm}.
We obtain then finally
$$a_{2n+2} +2 + \frac{q_{2n}-2}{q_{2n+1}}=p+1+\frac{\left((1-\beta')\beta^{n+1}-
    (1-\beta)\beta'^{n+1}\right)-2(\beta-\beta')}{\beta^{n+1}-\beta'^{n+1}}\,,$$
which is exactly $\text{ind}(w^{(n)})$.
This result holds for all parameters $p, q$ satisfying $p=q+1$, even for $p \leq 3$. Consequently, Theorem~\ref{ind} is in fact valid for all parameters $p,q$ with $p>q \geq 1$.

\section*{Acknowledgements}
The authors acknowledge financial support by Czech Science
Foundation GA \v{C}R 201/05/0169, by the grants MSM6840770039 and LC06002 of the
Ministry of Education, Youth, and Sports of the Czech Republic.


{\small

\end{document}